\documentclass[notitlepage,11pt]{article}
\usepackage{amssymb,amsmath,comment}
\catcode`\@=11
\@addtoreset{equation}{section}

\catcode`\@=12
\usepackage{colortbl}%

\def\R{\mathbb{R}}

\def\S{\mathbb{S}}
\def\f{\varphi}

\def\starstar{{2^{*\!*}}}
\def\thetastar{{\theta^{*\!*}}}

\def\irn{\int_{\R^n}}

\def\proof{\noindent{\textbf{Proof. }}}
\def\QED{\hfill {$\square$}\goodbreak \medskip}

\newtheorem{Theorem}{Theorem}[section]
\newtheorem{Lemma}[Theorem]{Lemma}
\newtheorem{Proposition}[Theorem]{Proposition}
\newtheorem{Corollary}[Theorem]{Corollary}
\newtheorem{Remark}[Theorem]{Remark}

\begin{document}

\title 
{Radially symmetric solutions to the\\  H\'enon-Lane-Emden system
on the critical hyperbola}

\author{Roberta Musina\footnote{Dipartimento di Matematica ed Informatica, Universit\`a di Udine,
via delle Scienze, 206 -- 33100 Udine, Italy. Email: {roberta.musina@uniud.it}}~ and
K. Sreenadh\footnote{Department of Mathematics, Indian Institute of Technology Delhi, Hauz Khas, New Delhi-110016, India.
Email:{sreenadh@maths.iitd.ac.in}}}

\date{}

\maketitle

\begin{abstract}
\footnotesize
We use variational methods to
study  the existence of nontrivial and radially symmetric solutions
to the H\'enon-Lane-Emden system with weights,
when the exponents involved lie on the "critical hyperbola". We
also discuss  qualitative properties of solutions
and nonexistence results.

\medskip

\medskip

\end{abstract}

\medskip

\section{Introduction}
\label{S:Introduction}

In this paper we discuss existence, non existence and qualitative properties of
nontrivial radially symmetric solutions $u,v$
to  the following weakly coupled system on the punctured
space $\R^n\setminus\{0\}$:
\begin{equation}
\label{eq:HLE}
\begin{cases}
-\Delta u=|x|^a|v|^{p-2}v\\
-\Delta v=|x|^{b}|u|^{q-2}u.
\end{cases}
\end{equation}
Here $n\ge 2$, $a, b\in\R$, $p,q>1$ belong to the
{\em weighted critical hyperbola}
\begin{equation}
\label{eq:hyperbola}
\frac{a+n}{p}+\frac{b+n}{q}=n-2,
\end{equation}
and satisfy the standard anticoercivity assumption
\begin{equation}
\label{eq:anticoercive}
\frac{1}{p}+\frac{1}{q}<1~\!.
\end{equation}
The H\'enon-Lane-Emden system
(\ref{eq:HLE}) is a largely studied problem. In the autonomous case $a=b=0$, P.L. Lions proved in  \cite{PLL}
the existence of a solution  $u\in\mathcal D^{2,q'}(\R^n), v\in \mathcal D^{2,p'}(\R^n)$
to
\begin{equation}
\label{eq:HLEaut}
\begin{cases}
-\Delta u=|v|^{p-2}v\\
-\Delta v=|u|^{q-2}u,
\end{cases}
\end{equation}
under the assumptions $n>n/p+n/q=n-2>0$. We quote also
the paper  \cite{HvdV} by Hulshof-Van der Vorst, for additional qualitative properties of the pair
$u,v$.

The role of the
{\em "critical hyperbola"} was first pointed out by Mitidieri
\cite{Mit0}, \cite{Mit1} for the autonomous case $a=b=0$ (see also \cite{SZ2}).
It turns out that (\ref{eq:HLEaut}) has no positive, radial solutions
$u,v\in C^2(\R^n)$
if $p,q$ are below the 
critical hyperbola. On the other hand, 
Serrin and Zou used   shooting methods
 in \cite{SZ}
to prove that 
(\ref{eq:HLE}) admits infinitely many positive radial solutions $u,v$ 
which tend to $0$ as 
$|x|\to +\infty$, provided that the pair $p,q$ is on or above the critical hyperbola. 

The {\em H\'enon-Lane-Emden conjecture} has been
raised in  \cite{SZ3} and in \cite{CDM} for 
a more general class of higher order system. It
says in particular that  there is no positive solution for system (\ref{eq:HLE})
if  $p,q$ is under the critical hyperbola. 
Bidaut-Veron and Giacomini have recently shown in \cite{BVG}
that if $n\ge 3$, $a,b>-2$, then the system (\ref{eq:HLE}) admits a positive 
classical radial solution $u,v$ with $u, v$
continuous at the origin
if and only if $(p, q)$ is above or on the critical hyperbola. 
We recall that
by Proposition 2.1 in  \cite{BV}, no solution to
(\ref{eq:HLE}) is continuous at the origin if $a\le-2$ or $b\le -2$.

Remarkable  results about the  H\'enon-Lane-Emden conjecture have been recently obtained 
also in  \cite{BusMan}, \cite{F},
\cite{FG}, \cite{Ph}, \cite{PQS}, \cite{SZ2} and \cite{Sou}. 

Finally,
we recall that the weighted critical hyperbola
enters in a natural way in the context of the solvability of Hardy-H\'enon 
type elliptic systems in bounded domains, see for
instance \cite{dFPR08} and \cite{CR}.

In the present paper we first use variational methods to extend the Lions existence theorem
to the non autonomous case. Then we discuss
 nonexistence results. We always assume that (\ref{eq:hyperbola}) and (\ref{eq:anticoercive})
are satisfied.
We limit ourselves to state here some of our results, and we postpone
more precise statements to Section \ref{S:HLE}.

\medskip
\noindent
{\bf Existence} If 
$a\neq-n$ and $b\neq -n$, then (\ref{eq:HLE}) has a nontrivial radial solution $u,v$
such that
$$
\int_{\R^n}|x|^{-\frac{a}{p-1}}|\Delta u|^{p'}~\!dx<\infty~,\quad \int_{\R^n}|x|^{-\frac{b}{q-1}}|\Delta v|^{q'}~\!dx<\infty~\!.
$$
Moreover, it holds that
\begin{equation}
\label{eq:0}
\lim_{|x|\to \infty}|x|^{\frac{b+n}{q}} u(x)=
\lim_{x\to 0}|x|^{\frac{b+n}{q}} u(x)=
\lim_{|x|\to \infty}|x|^{\frac{a+n}{p}}v(x)=
\lim_{x\to 0}|x|^{\frac{a+n}{p}}v(x)=0,
\end{equation}
and $u$, $v$ are both positive if and only if $a>-n$ and $b>-n$.

 \medskip
\noindent
 {\bf Nonexistence}
 Let $u,v$ be a solution to (\ref{eq:HLE}) on $\R^n\setminus\{0\}$ and assume
that either
\begin{equation}
\label{eq:uvenergy}
\begin{cases}
~\displaystyle\lim_{|x|\to \infty}|x|^{\frac{b+n}{q}} u(x),~ \displaystyle\lim_{|x|\to \infty}|x|^{\frac{a+n}{p}}v(x)
&\textit{exist and are finite, or}\\
~~\displaystyle\lim_{x\to 0}~|x|^{\frac{b+n}{q}} ~\!u(x),~
~\displaystyle\lim_{x\to 0}~|x|^{\frac{a+n}{p}}~\!v(x)
&\textit{exist and are finite.}
\end{cases}
\end{equation}
If $a\le-n$ or $b\le-n$ and if $u\ge 0, v\ge 0$ then $u\equiv v \equiv 0$.

\medskip

Let us briefly describe our approach. It has already been noticed for instance
in \cite{BV}, \cite{BVG}, \cite{BusMan}, that radial solutions to (\ref{eq:HLE}) are in one-to-one correspondence
with trajectories $g,f$ of the Hamiltonian system
\begin{equation}
\label{eq:ODEsystem}
\begin{cases}
-g''+2Ag'+\Gamma g=|f|^{p-2}f&\textrm{on $\R$}\\
-f''-2Af'+\Gamma f=|g|^{q-2}g&\textrm{on $\R$}
\end{cases}
\end{equation}
for suitable constants $A,\Gamma\in\R$ depending on the data.
Notice that (\ref{eq:ODEsystem}) includes the Schr\"odinger equation
$-g''+\Gamma g=|g|^{p-2}g$,
whose relevance with the Caffarelli-Kohn-Nirenberg inequality
was pointed out by Catrina and Wang in \cite{CatWan}. For $p=2$ the system (\ref{eq:ODEsystem})
reduces to the fourth order ordinary differential equation
\begin{equation}
\label{eq:ODEeq}
g''''-2(2A^2+\Gamma)g''+\Gamma ^2 g=|g|^{q-2}g~\!,
\end{equation}
which is naturally related to
second order dilation invariant inequalities of Rellich-Sobolev type, see \cite{BM}. 
Actually the system (\ref{eq:ODEsystem}) and the equation (\ref{eq:ODEeq}) have  
independent interest because of their applications. We shall not attempt to give 
a complete list of references. We cite for instance \cite{Bu}
 \cite{CT}, \cite{CMcK}
\cite{GG}, 
\cite{H1}--\cite{H3}, \cite{LS}  and references therein. In the monograph \cite{PT} by Peletier and Troy one can find
several applications and a rich bibliography on these topics.

In Section \ref{S:ODE} we use the results in \cite{Mu12} and variational methods
to get the existence of solutions $g\in W^{2,p'}(\R), f\in W^{2,q'}(\R)$ to (\ref{eq:ODEsystem}), compare
with Theorem \ref{T:1dim}. 
Then we discuss sign properties of  solutions to (\ref{eq:ODEsystem}) having certain
behavior at $-\infty$ and/or at $+\infty$. In Section \ref{S:HLE} we obtain
our main theorems about (\ref{eq:HLE}) as corollaries of our results for (\ref{eq:ODEsystem}). 

In the Appendix we indicate a possible
non-radial approach to (\ref{eq:HLE}).

\bigskip

\small
\noindent
{\bf Notation} For any integer $n\ge 2$ we denote by
$\omega_n$ the $n-1$ dimensional measure of the unit sphere
$\S^{n-1}$ in $\R^n$.

Let $q\in[1,+\infty)$ and let $\omega$ be a non-negative measurable function on
a domain $\Omega\subseteq\R^n$, $n\ge 1$. The weighted
Lebesgue space $L^q(\Omega;\omega(x)~\! dx)$ is the space of measurable maps
$u$ in $\Omega$ with finite norm $\left(\int_{\Omega}|u|^q
\omega(x)~\!dx\right)^{1/q}$.  For $\omega\equiv 1$ we simply write
$L^q(\Omega)$. As usual, $\|\cdot\|_\infty$ is the $L^\infty$ norm.

For any function $\f:\R\to \R$, the notation $\f(\pm\infty)=c$ means that there exists
$\displaystyle\lim_{s\to\pm\infty}\f(s)=c$.

\normalsize

\section{A $2\times2$ system of ordinary differential equations} 
\label{S:ODE}

In this section we provide conditions for the existence of  solutions to (\ref{eq:ODEsystem}) vanishing at $\pm\infty$
and for the nonexistence of positive solutions
having nonnegative limits at $-\infty$ or at $\infty$.
We start with an existence result.
\begin{Theorem}
\label{T:1dim}
Let $p, q\in(1,\infty)$ $A,\Gamma\in\R$ be given,
such that $A^2+\Gamma\ge 0$ and $\Gamma\neq 0$.
Assume that  (\ref{eq:anticoercive}) is satisfied.
Then the system (\ref{eq:ODEsystem}) has a nontrivial
solution $g, f$ such that
$g\in W^{2,p'}(\R)$ and $f\in W^{2,q'}(\R)$.
\end{Theorem}

\proof
To simplify notations, we set
$$
\mathcal L_+\f:= -\f''+2A\f'+\Gamma \f~,\quad \mathcal L_-\f=-\f''-2A\f'+\Gamma \f.
$$
Since $A^2+\Gamma \ge 0$ and $\Gamma \neq 0$, from Proposition
5.2 in \cite{Mu12} we have that the infimum
\[
I_{p',q}(A,\Gamma)=
\inf_{\scriptstyle g\in W^{2,p'}(\R)\atop\scriptstyle g\ne 0}
\frac{\displaystyle \int_{\R}|\mathcal L_+g|^{p'} ds}
{\left(\displaystyle\int_{\R}|g|^q dx\right)^{p'/q}}
\]
is achieved by some
$g\in W^{2,p'}(\R)$ that solves
$$
\int_\R|\mathcal L_+g|^{p'-2}\mathcal L_+g\mathcal L_+\psi~ds=\int_\R|g|^{q-2}g\psi~ds\quad
\textrm{for any $\psi\in W^{2,p'}(\R)$.}
$$
Thus $g\in W^{2,p'}(\R)$ is a weak solution to the following fourth-order ODE:
$$
\mathcal L_-\left(|\mathcal L_+g|^{p'-2}\mathcal L_+g\right)=|g|^{q-2}g\quad
\textrm{on $\R$,}
$$
which is equivalent to the system (\ref{eq:ODEsystem}), by defining
$f=-|\mathcal L_+g|^{p'-2}\mathcal L_+g$.
Clearly, $g,f\in C^2(\R)$. Now we recall that
$\|\mathcal L_-\cdot\|_{q'}$ is an equivalent norm in $W^{2,q'}(\R)$ by
Proposition 5.2 in \cite{Mu12}. Since $g\in W^{2,p'}(\R)\hookrightarrow L^q(\R)$, we have
$\mathcal L_-f=|g|^{q-2}g\in L^{q'}(\R)$, and thus $f\in W^{2,q'}(\R)$.
\QED

\begin{Remark}
One could exchange $g$ and $f$ in the proof
of Theorem \ref{T:1dim} to find a solution $\tilde g,\tilde f$, such that
$\tilde f\in W^{2,q'}(\R)$ achieves $I_{q',p}(-A,\Gamma)$ and
$\tilde g=-|\mathcal L_-\tilde f|^{p'-2}\mathcal L_-\tilde f$. This
argument does not lead to a multiplicity result for (\ref{eq:ODEsystem}).
To simplify notations we set
$m=I_{p',q}(A,\Gamma)$ and $\tilde m=I_{q',p}(-A,\Gamma)$. Since
$|\mathcal L_-f|^{q'}=|g|^{q}$, $|f|^p=|\mathcal L_+g|^{p^\prime}$, and since $g$ achieves $m$ we find
$$
\tilde m\le\frac{\displaystyle\int_\R|\mathcal L_-f|^{q'}~\!ds}
{\left(\displaystyle\int_\R|f|^p~\!ds\right)^{q'/p}}=
\frac{\displaystyle\int_\R|g|^q~\!ds}
{\left(\displaystyle\int_\R|\mathcal L_+g|^p~\!ds\right)^{q'/p}}=
m^{\frac{p-q'}{p}~\!\frac{q}{q-p'}},
$$
so that $\tilde m^{\frac{q-p'}{q}}\le m^{\frac{p-q'}{p}}$.
In a similar way we get the opposite inequality, and in particular
$\tilde m^{\frac{q-p'}{q}}= m^{\frac{p-q'}{p}}$. Moreover, $\tilde f$ achieves
$\tilde m$ and $\tilde g$ achieves $m$.
\end{Remark}

In order to study the qualitative properties of solutions to (\ref{eq:ODEsystem}) we
 take advantage of its Hamiltonian structure. Indeed,
the system (\ref{eq:ODEsystem}) is conservative, and
any  solution $g, f$ satisfies
\begin{equation}
\label{eq:conservation}
E(g,f):=g'f'-\Gamma gf+\frac{1}{q}|g|^q+\frac{1}{p}|f|^p= \textrm{constant}.
\end{equation}

\begin{Remark}
\label{R:regularity}
Let $g\in W^{2,p'}(\R)$, $f\in W^{2,q'}(\R)$ be a solution to
(\ref{eq:ODEsystem}). By well known facts about Sobolev spaces, the functions
$g, g', f$ and $f'$ are H\"older continuous on $\R$.  Thus
$g,f\in C^2(\R)$.  In addition $g, g', f$ and $f'$ vanish
at $\pm\infty$ and hence  (\ref{eq:conservation}) implies
\begin{equation}
\label{eq:conservation_null}
g'f'-\Gamma gf+\frac{1}{q}|g|^q+\frac{1}{p}|f|^p\equiv 0\quad \textrm{on $\R$.}
\end{equation}
\end{Remark}

\begin{Remark}
\label{R:Ham}
Problem (\ref{eq:ODEsystem}) is equivalent to a $(2\times2)$-dimensional
first order Hamiltonian system. For $X=(x_1,x_2), Y=(y_1,y_2)\in\R^2$
we set
$$
H(X,Y)=y_1y_2+A(x_1y_1-x_2y_2)-(A^2+\Gamma)x_1x_2+
\frac{1}{q}|x_1|^q+\frac{1}{p}|x_2|^p.
$$
Then $g, f$ solves (\ref{eq:ODEsystem}) if and only if
$X=(g,f), Y=(f'+Af,g'-Ag)$ solves
\begin{equation}
\label{eq:Ham}
\begin{cases}
X'=\partial_YH(X,Y)\\
Y'=-\partial_XH(X,Y).
\end{cases}
\end{equation}
If $\Gamma\neq 0$ and $\delta:={pq-(p+q)}>0$,
then $
\pm~\!\left(|\Gamma|^{p/\delta},|\Gamma|^{-1+q/\delta}\Gamma\right)
$
are equilibrium points for (\ref{eq:Ham}). Notice that a positive equilibrium exists
if and only if $\Gamma>0$.
\end{Remark}

From (\ref{eq:conservation}) we first infer the following a-priori bound on trajectories having
null energy.

\begin{Proposition}
\label{P:aprioribound}
Let $g,f\in C^2(\R)$ be a solution to (\ref{eq:ODEsystem}) such that
$g, g', f$ and $f'$ vanish at infinity. Then
$$
\|g\|_\infty^{q-p'}\le \frac{q}{p'}|\Gamma|^{p'}~,\quad
\|f\|_\infty^{p-q'}\le \frac{p}{q'}|\Gamma|^{q'}.
$$
In particular, if $\Gamma=0$ then $g=f\equiv 0$.
\end{Proposition}

\proof
Let $\overline s\in\R$ be such that $|g(\overline s)|=\|g\|_\infty$. Then
$g'(\overline s)=0$ and therefore from (\ref{eq:conservation}) and since $E(g,f)=0$ we get
$$
\frac{1}{q}\|g\|_\infty^q+\frac{1}{p}|f(\overline s)|^p=
\Gamma f(\overline s)\|g\|_\infty\le
\frac{|\Gamma|^{p'}}{p'}\|g\|_\infty^{p'}+\frac{1}{p}|f(\overline s)|^p
$$
by Young's inequality. The desired a-priori bound on $g$ follows immediately.
The  estimate on $\|f\|_\infty$ can be obtained in a similar way.
\QED
In the remaining part of this section,  we study the
sign of solutions $g,f$ to (\ref{eq:ODEsystem}). We distinguish the
case $\Gamma>0$ from the case when $\Gamma$ is nonpositive.

\begin{Theorem}
\label{T:change_sign}
Let $g, f\in C^2(\R)$ be a solution to (\ref{eq:ODEsystem}), such
that $g$ and $f$ vanish at $\pm \infty$ together with their derivatives.
If $\Gamma>0$
then $g\equiv f\equiv 0$ or $gf>0$ on $\R$.
\end{Theorem}

\proof
We start by noticing that $g,f$ satisfies (\ref{eq:conservation_null}).
In a moment we will prove the following
\bigskip

\noindent
{\bf Claim:} {\em $g(s)f(s)\neq 0$ for any $s\in\R$.}

\medskip

Assume  that the claim is proved.
Then both $g$ and $f$ have constant sign. The function $g$ has at least
one critical point $\overline s$. By (\ref{eq:conservation_null}), it holds that
$$
-\Gamma g(\overline s)f(\overline s)+\frac{1}{q}|g(\overline s)|^q
+\frac{1}{p}|f(\overline s)|^p=0.
$$
Thus $g(\overline s)f(\overline s)>0$, and therefore $gf>0$ everywhere in $\R$, that
concludes the proof of the theorem.

\medskip

It remains to prove the claim. Notice that from
(\ref{eq:conservation_null}) the following facts follow:
\begin{gather}
\label{eq:fact}
\textit{if ~$g'(\xi)f'(\xi)=0$~ then ~ $f(\xi)=g(\xi)=0$~ or~ $f(\xi)g(\xi)>0$}\\
\label{eq:factnew}
\textit{if ~$g(\xi)f(\xi)=0$~~then ~
$f(\xi)=g(\xi)=0=f'(\xi)g'(\xi)$~ or ~ $f^\prime(\xi)g^\prime (\xi)<0$.}
\end{gather}
By contradiction, assume that $g$
vanishes somewhere. Up to a change of sign and/or inversion
$s\mapsto -s$ we can assume that $g$ attains its negative minimum
at some $s_1\in \R$ and that $g$ reaches $0$ in $(s_1,\infty)$.
Let $s_2$ be the first zero of $g$ in $(s_1,\infty)$. Thus
$g<0$ on $[s_1,s_2)$,  $f(s_1)<0$ by
(\ref{eq:fact}), and $g'(s_2)\ge 0$. In addition,
\begin{equation}
\label{eq:fact3}
\textit{if $f'(\bar s)=0$~ for some~
$\bar s\in [s_1,s_2)$,  then~ $f(\bar s)<0$,}
\end{equation}
because of (\ref{eq:fact}). Now we prove that
\begin{equation}
\label{eq:f_negative}
g'(s_2)f'(s_2)=0~,\quad f(s_2)=0~,\quad f<0\quad\textit{on $[s_1,s_2)$.}
\end{equation}
If $g'(s_2)=0$ then (\ref{eq:f_negative}) readily follows from (\ref{eq:fact}) and (\ref{eq:fact3}).
If $g'(s_2)>0$ and $f'(s_2)=f(s_2)=0$ then (\ref{eq:fact3}) immediately implies
(\ref{eq:f_negative}).
In view of (\ref{eq:factnew}), to conclude the proof of (\ref{eq:f_negative})
we only have to exclude that $g'(s_2)>0>f'(s_2)$. We argue by contradiction.
If $f'(s_2)<0$ then $f(s_2)<0$ by (\ref{eq:fact3}). Since $g$ is increasing
in a neighborhood of $s_2$ and since $g$ decays at infinity, there is a point
$s_3>s_2$ such that $g'(s_3)=0$ and $g>0$ on $(s_2,s_3]$. But then
$f(s_3)>0$ by (\ref{eq:fact}). Since $f(s_2), f'(s_2)$ are negative, we infer that $f$ has a minimum
$s_4\in(s_2,s_3)$, with $f(s_4)<0$. But then $g(s_4)<0$ by (\ref{eq:fact}),
which is impossible. Thus (\ref{eq:f_negative}) is proved.

In conclusion, we have that $g, f$ solves the system
$$
\begin{cases}
g''-2Ag'-\Gamma g=-|f|^{p-2}f\ge 0~&\textrm{in $(s_1,s_2)$}\\
f''+2Af'-\Gamma f=-|g|^{q-2}g\ge 0&\textrm{in $(s_1,s_2)$}\\
g, f<0&\textrm{in $(s_1,s_2)$}\\
g(s_2)=f(s_2)=g'(s_2)=f'(s_2)=0~\!,
\end{cases}
$$
that contradicts the Hopf boundary
point lemma.
The claim and  the theorem are completely proved.
\QED

The condition $\Gamma>0$ is also necessary to have the
existence of positive  solutions vanishing at $\pm\infty$. In view of Remark
\ref{R:regularity}, the next proposition applies in particular to solutions
$g\in W^{2,p'}(\R), f\in W^{2,q'}(\R)$.

\begin{Proposition}
\label{P:new}
Let $g, f\in C^2(\R)$ be a solution to (\ref{eq:ODEsystem}), such
that $g$ and $f$ vanish at $\pm \infty$ together with their derivatives.
If $\Gamma\le 0$ and $gf\ge 0$ on $\R$ then
$g\equiv f\equiv 0$.
\end{Proposition}

\proof
The  trajectory $g,f$ has null energy, that is,
(\ref{eq:conservation_null}) holds. In particular, at any critical point
$\bar s$ of $g$ one has that
$|\Gamma| g(\bar s)f(\bar s)+\frac{1}{q}|g(\bar s)|^q+\frac{1}{p}|f(\bar s)|^p=0$.
Thus both $g$ and $f$ vanish at $\bar s$. In particular, $\min g=\max g=0$,
and the conclusion follows.
\QED
 We conclude this section with two more nonexistence results in case
 $\Gamma\le 0$.

\begin{Theorem}
\label{T:perBVG}
Assume that $g,f\in C^2(\R)$ solves (\ref{eq:ODEsystem}) for some
$A\in \R$, $\Gamma\le 0$ and $p,q$ satisfying (\ref{eq:anticoercive}). In addition, assume that
$$
g(-\infty)= c_g\in [0,\infty)~,\quad
f(-\infty)= c_f\in [0,\infty)~,\quad g\ge 0~\textit{and}~f\ge 0\quad\textit{on $\R$.}
$$
Then $g\equiv f\equiv 0$.
\end{Theorem}

\proof
First of all we notice that $g,f$ can not be a nontrivial pair of constant functions
by Remark \ref{R:Ham}.

The function $h:=-f'-2Af$ is increasing
in $\R$, as $h'=g(s)^{q-1}-\Gamma f\ge 0$. Thus it has a limit as $s\to -\infty$.
Hence, also $f'$ has a limit as $s\to -\infty$. Clearly
\begin{equation}
\label{eq:limitf'}
f'(-\infty)=0,
\end{equation}
and therefore from (\ref{eq:ODEsystem}) we get also
\begin{equation}
\label{eq:limitf''}
-f''(-\infty)=-\Gamma c_f+c_g^{q-1}\ge 0.
\end{equation}
In a similar way we get that
\begin{equation}
\label{eq:limitg}
g'(-\infty)=0~,\quad-g''(\infty)=-\Gamma c_g+c_f^{p-1}\ge 0.
\end{equation}
In particular, from (\ref{eq:conservation}) and (\ref{eq:limitf'}), (\ref{eq:limitg}) we infer that
$$
g'f'-\Gamma gf+\frac{1}{q}|g|^q+\frac{1}{p}|f|^p=-\Gamma c_g c_f+\frac{1}{q} c_g^q+
\frac{1}{p} c_f^p\quad\textrm{on $\R$}.
$$

\medskip

\noindent
{\bf Claim 1:} {\em If $ c_g= c_f=0$ then $g\equiv f\equiv 0$.}

\smallskip
\noindent
To prove the claim, we notice that
the trajectory $g,f$ satisfies (\ref{eq:conservation_null}). 
If we assume by contradiction that $g$ or $f$ do not vanish identically,
then there exists $s_0\in\R$ such that $g'(s_0)f'(s_0)<0$. To fix ideas,
assume that $f'(s_0)<0$. Since $f\ge 0$ and $f(s)\to 0$ as $s\to -\infty$,
it means that $f$ must have a positive local maximum $s_1<s_0$. At the point
$s_1$ the conservation law (\ref{eq:conservation_null}) gives
$-\Gamma g(s_1)f(s_1)+\frac{1}{q}|g(s_1)|^q+\frac{1}{p}|f(s_1)|^p=0$,
which contradicts $f(s_1)>0$. The claim is proved.
\bigskip

\noindent
{\bf Claim 2} {\em If $A\le 0$ then  $\Gamma c_f=0~~\textrm{and}~~ c_g=0$}

\medskip
\noindent
By contradiction, assume that $-\Gamma c_f+ c_g^{q-1}>0$.
Then the function $f$ is strictly concave and
decreasing in a neighborhood of $-\infty$ by (\ref{eq:limitf''}) and (\ref{eq:limitf'}). Thus in particular $ c_f>0$, and therefore form the conservation law we get
\begin{equation}
\label{eq:contradiction1}
g'f'-\Gamma gf+\frac{1}{q}|g|^q+\frac{1}{p}|f|^p\ge \frac{1}{p} c_f^p>0
\quad\textrm{on $\R$}.
\end{equation}
Since $f$ is bounded from below, it can not be strictly concave on $\R$.
We claim that $f$ can never  be locally convex.
Assume that there exists $s_0\in\R$ such that $f''(s_0)>0$.
Then from (\ref{eq:ODEsystem})
we have that
$-2A f'(s_0)>-\Gamma f(s_0)+g(s_0)^{q-1}\ge 0$.
Thus, $A<0$ and $f'(s_0)>0$. Since $f'(s)<0$ for $s<<0$, then the function $f$ must have
a local minimum $s_1\in(-\infty, s_0)$. Thus $f'(s_1)=0$ and $f''(s_1)\ge 0$. But then
$$
0\ge -f''(s_1)=-\Gamma f(s_1)+g(s_1)^{q-1}\ge 0,
$$
which implies $\Gamma f(s_1)=g(s_1)=0$. In particular, $g'(s_1)=0$, and $g''(s_1)\ge 0$,
since $s_1$ is a minimum for $g$ thanks to the assumption that $g\ge 0$. Thus, (\ref{eq:ODEsystem}) gives
$0\ge -g''(s_1)=f(s_1)^{p-1}\ge 0$. Thus
 $f(s_1)=0$, contradicting (\ref{eq:contradiction1}).

We have proved that $f''\le 0$ on $\R$.
Thus there exists $s_0\in\R$
such that $f$ is a nonnegative constant on $[s_0,\infty)$.
But then from (\ref{eq:ODEsystem}) we infer that
$f\equiv g\equiv 0$ on $[s_0,\infty)$, as $\Gamma\le 0$.
We have reached again a contradiction with (\ref{eq:contradiction1}),
and the claim is proved.

\bigskip

\noindent
{\bf Claim 3} {\em If $A\ge 0$ then  $\Gamma c_g=0~~\textrm{and}~~ c_f=0$}

\medskip
\noindent
It is sufficient to exchange the roles of $g$ and $f$, and argue as in Claim 2.

\bigskip

Now we are in position to conclude the proof. By Claim 1, we only have to show that 
$ c_g= c_f=0$. Thus we are done
if $A=0$, thanks to Claims 2 and 3. We have to study the case
\begin{equation}
\label{eq:last_contradiction}
A<0~,\quad\Gamma= c_g=0
\end{equation}
and the case $A>0$, $\Gamma=0= c_f$, that can can be handled in a similar way.
Assume that (\ref{eq:last_contradiction}) holds. Since $g$ solves
$-g''+2Ag'=f^{p-1}\ge 0$,
then the function $-g'+2Ag$ is non decreasing on $\R$.
Hence $-g'+2Ag\ge 0$  by (\ref{eq:limitg}) and since $ c_g=0$. Thus $g'\le 2Ag\le 0$ on $\R$
that is, $g\equiv 0$ because it is non increasing and nonnegative.
The proof is complete.
\QED

Since the system (\ref{eq:ODEsystem}) is invariant with respect to
inversion $s\mapsto -s$, then clearly the next result holds as well.

\begin{Theorem}
\label{T:perBVG_symm}
Assume that $g,f\in C^2(\R)$ solves (\ref{eq:ODEsystem}) for some
$A\in\R$, $\Gamma\le 0$ and $p,q$ satisfying (\ref{eq:anticoercive}). In addition, assume that
$$
g(\infty)= c_g\in [0,\infty)~,\quad
f(\infty)= c_f\in [0,\infty)~,\quad g\ge 0~\textit{and}~f\ge 0\quad\textit{on $\R$.}
$$
Then $g\equiv f\equiv 0$.
\end{Theorem}

We conclude this section with a result that holds in case
$p=2<q$.

\begin{Theorem}
\label{T:BM}
Let $q\in(2,\infty)$ and assume that
$A^2+\Gamma \ge 0~,\quad \Gamma \neq 0$. Up to translations in $\R$,
composition with the inversion $s\mapsto -s$ and change of sign, the system
$$
\begin{cases}
-g''+2Ag'+\Gamma g=f&\textrm{on $\R$}\\
-f''-2Af'+\Gamma f=|g|^{q-2}g&\textrm{on $\R$.}
\end{cases}
$$
has a unique nontrivial solution $(g,f)$ such that $g\in H^{2}(\R)$
and $f\in W^{2,q'}(\R)$.
Moreover,  $g$ is even, positive and strictly decreasing on $(0,\infty)$,
and $f$ is positive if and anly if $\Gamma>0$.
\end{Theorem}

\proof
Existence is given by Theorem \ref{T:1dim}.
Notice that $g$ is smooth and solves
\begin{equation}
\label{eq:fourthODE}
g''''-2(2A^2+\Gamma)g''+\Gamma ^2 g=|g|^{q-2}g.
\end{equation}
On the other hand, since $(2A^2+\Gamma)^2\ge \Gamma ^2$, then Theorem 2.2 in \cite{BM}
implies that (\ref{eq:fourthODE}) has a unique solution $g$ (up to the above
transforms), that can be taken to be positive, even and strictly decreasing on $(0,\infty)$.
The uniqueness of $f$ is immediate. The last statement concerning the sign of
$f$ follows by Theorem \ref{T:change_sign} and Proposiion \ref{P:new}.
\QED

\begin{Remark}
Clearly,
$f$ is even if and only if $A$=0.
\end{Remark}

\section{The H\'enon-Lane-Emden system}
\label{S:HLE}
In this section we
provide conditions for the existence of solutions to
\eqref{eq:HLE}
in suitable energy spaces
and for the nonexistence of positive solutions
having certain behavior at $0$ or at $\infty$.

We start by introducing some weighted Sobolev spaces.
Let $\theta\in(1,\infty)$ and $\alpha\in\R$ be given, such that
$\alpha\notin\{ 2\theta-n, np-n\}$.
Then we can use the results in \cite{Mu12} to define the Banach space
$\mathcal D^{2,\theta}_{\rm r}(\R^n;|x|^\alpha dx)$ as the completion of radial functions
in $C^2_{c}(\R^n\!\setminus\!\{0\})$ with respect to the norm
$$
\|u\|_{\alpha}=\left(\irn|x|^\alpha|\Delta u|^\theta~dx\right)^{1/\theta}.
$$
To any pair of radial
functions $u,v\in C^2(\R^n\setminus\{0\})$, we associate the pair
$g,f\in C^2(\R)$ defined by
\begin{equation}
\label{eq:EFtransform}
u(x)=|x|^{-\lambda_1}~\!g\left(-\log|x|\right)~,\quad
v(x)=|x|^{-\lambda_2}~\!f\left(-\log|x|\right)~\!,
\end{equation}
where
$$
\lambda_1=\frac{b+n}{q}~,\quad
\lambda_2=\frac{a+n}{p}~\!.
$$
We will always assume that $(p,q)$ belongs to the {\em critical hyperbola} in (\ref{eq:hyperbola}),
that is, $\lambda_1+\lambda_2=n-2$.

We introduce also the constants
\begin{equation}
\label{eq:defAGamma}
\Gamma=\frac{n+a}{p}~\!\frac{n+b}{q}=\lambda_1\lambda_2~,
\quad
A=\frac{n-2}{2}-\lambda_1=-~\!\frac{n-2}{2}+\lambda_2~\!.
\end{equation}
Notice that
$$
A^2+\Gamma=\left(\frac{n-2}{2}\right)^2\ge 0.
$$
A direct computation shows that a radial pair
$u,v$ solves  (\ref{eq:HLE}) on $\R^n\setminus\{0\}$ if and only if $g,f$ solves
(\ref{eq:ODEsystem}) with $\Gamma$, $A$ given by (\ref{eq:defAGamma}).
Thanks to the results in previous section we first get the next existence theorem.

\begin{Theorem}
\label{T:HLE}
Let $n\ge 2$, $a,b\in\R\setminus\{-n\}$ and $p,q>1$. Assume that (\ref{eq:anticoercive})
and (\ref{eq:hyperbola}) are satisfied. Then the H\'enon-Lane-Emden system (\ref{eq:HLE}) has a
radially symmetric solution
\begin{equation}
\label{eq:uv_space}
u\in \mathcal D^{2,p'}_{\rm r}(\R^n;|x|^{-\frac{a}{p-1}} dx)~,\quad
v\in \mathcal D^{2,q'}_{\rm r}(\R^n;|x|^{-\frac{b}{q-1}} dx).
\end{equation}
 In addition, $u, v$ satisfies (\ref{eq:0})
\end{Theorem}

\proof
Define  $\Gamma$, $A$ as in (\ref{eq:defAGamma}), and notice that $\Gamma\neq 0$, and $A^2+\Gamma\ge 0$.
By Theorem \ref{T:1dim}, we see that there exists $f \in W^{2,p'}(\R)$ and $g \in W^{2,q'}(\R)$ satisfying \eqref{eq:ODEsystem}.
Now using the Emden-Fowler transformation in \eqref{eq:EFtransform} and the results in \cite{Mu12}, we get $u, v$ satisfies (\ref{eq:uv_space})
and solves (\ref{eq:HLE}). The conclusion readily follows since
$$
|x|^\frac{b+n}{q} u(x)=g(-\log|x|)~,\quad
|x|^\frac{a+n}{p} v(x)=f(-\log|x|)
$$
and since $g, f$ vanish at $\pm \infty$.
\QED

\begin{Theorem}
\label{T:ne}
Let $n\ge 2$, $a,b\in\R$ and $p,q>1$. Assume that (\ref{eq:anticoercive})
and (\ref{eq:hyperbola}) are satisfied.
Let $u,v\in C^2(\R^n\setminus\{0\})$ be a radially symmetric solution to (\ref{eq:HLE})
on $\R^n\setminus\{0\}$.
\begin{description}
\item$i)$
If $a>-n$, $b>-n$ and if $u,v$  satisfies (\ref{eq:0}), then $u\equiv v \equiv 0$ or $uv>0$ on $\R$.

\item$ii)$
Assume that (\ref{eq:uvenergy}) holds.
If $a\le-n$ or $b\le-n$ and if $u\ge 0, v\ge 0$ then $u\equiv v \equiv 0$.

 \end{description}
 \end{Theorem}

\proof
Define $A,\Gamma$ and use the Emden-Fowler transform
$(u,v)\mapsto (g,f)$ as before.
Notice that $\Gamma>0$ in case $i)$ and
$\Gamma\le 0$ in case $ii)$. Then apply Theorem \ref{T:perBVG}
and Theorem \ref{T:change_sign}.
\QED

In the next corollary we emphasize the impact of Theorem \ref{T:ne} in case $n=2$,
when Theorem \ref{T:HLE} gives existence on the critical hyperbola whenever
$a, b\neq -2$.

\begin{Corollary}
\label{C:n=2}
Let $n= 2$ and $p,q>1$. Assume that (\ref{eq:anticoercive})
and (\ref{eq:hyperbola}) are satisfied, and in addition assume that
$a, b\neq -2$. Let $u,v\in C^2(\R^2\setminus\{0\})$ be a radially
symmetric and nonnegative solution to (\ref{eq:HLE}) satisfying
(\ref{eq:uvenergy}). Then $u\equiv v \equiv 0$.
\end{Corollary}

In Theorem \ref{T:ne} we saw that the sign of $\Gamma$ affects the sign of
the product $uv$. However, at least in case $p=2$, the function $u$
never changes sign, also in case $\Gamma<0$. The next result for problem
\begin{equation}
\label{eq:p=2}
\begin{cases}
-\Delta u=|x|^a v\\
-\Delta v=|x|^{b}|u|^{q-2}u.
\end{cases}
\end{equation}
is an immediate
consequence to Theorem \ref{T:BM}.

\begin{Theorem}
\label{T:CorBM}
Let $n\ge 2$, $a,b\in\R$ and $q>1$. Assume that $a, b\neq -n$ and
$$
\frac{a+n}{2}+\frac{b+n}{q}=n-2.
$$
is satisfied. Up to dilations, compositions with the Kelvin
transform and change of sign,  problem (\ref{eq:p=2})
has a unique nontrivial radial solution
$u\in \mathcal D^{2,2}_{\rm r}(\R^n;|x|^{-a} dx)$, $v\in \mathcal D^{2,q'}_{\rm r}(\R^n;|x|^{-\frac{b}{q-1}} dx)$.
Moreover,
$u$ is positive, and $v$ is positive if and only if $a, b>-n$.
\end{Theorem}

\appendix

\noindent

\section{\!\!\!\!\!\!ppendix: a non-radial approach}
\label{A:nonradial}

Following
Wang \cite{Wa} and Calanchi-Ruf   \cite{CR}, we notice that  (\ref{eq:HLE})
is formally equivalent to the  fourth order equation
\begin{equation}
\label{eq:equation}
\Delta(|x|^\alpha|\Delta u|^{\theta-2}\Delta u)=|x|^{b}|u|^{q-2}u
\end{equation}
where
$\theta=p'=\frac{p}{p-1}$ and $\alpha=-\frac{a}{p-1}$. 
Equation (\ref{eq:equation}) is variational. In particular, its nontrivial solutions
can be find as critical points for the functional
$$
u~\mapsto~\frac{\displaystyle
\int_{\R^n}|x|^{\alpha}|\Delta u|^\theta dx}
{\displaystyle\left(\int_{\R^n}|x|^{b}|u|^{q}dx\right)^{\theta /q}}~\!
$$
on a suitable function space. Let us introduce the weighted Rellich constant
\begin{equation}
\label{eq:Mit_ine}
\mu_{\theta }(\alpha):  =  \inf_{\scriptstyle u\in C^{2}_{c}(
\R^n\setminus\{0\})
\atop\scriptstyle u=u(|x|)~,~ u\ne 0}\frac{\displaystyle
\irn |x|^{\alpha}|\Delta u|^\theta dx}
{\displaystyle\irn |x|^{\alpha-2\theta}|u|^{\theta }dx}~\!.
\end{equation}
The  best constant $\mu_{\theta }(\alpha)$ is explicitly known in few cases. We define
$\Gamma$ as in (\ref{eq:defAGamma}) and we notice that
\begin{equation}
\label{eq:gamma}
\Gamma =
\left(\frac{n+\alpha}{\theta}-2\right)\left(n-\frac{n+\alpha}{\theta}\right)
\end{equation}
if  (\ref{eq:hyperbola}) is satisfied.
The value of $\mu_{2 }(\alpha)$ (case $\theta=2$) is known from \cite{GM}, \cite{CM1}: 
$$
\mu_{2 }(\alpha)=\min_{k\in\mathbb N\cup\{0\}}\left|\Gamma+k(n-2+k)\right|^2.
$$
For general $\theta>1$, Mitidieri proved in \cite{Mit2} that 
$\mu_{\theta }(\alpha)=\left|\Gamma\right|^\theta$, provided that $\Gamma\ge0$.

From now on, we assume that
$\mu_{\theta }(\alpha)>0$. Then we can define the space
$\mathcal D^{2,\theta}(\R^n;|x|^\alpha dx)$ as the closure of functions in $C^{2}_{c}(\R^n\backslash\{0\})$ 
with respect to the norm
\[\|u\|^\theta =\int_{\R^n} |x|^\alpha |\Delta u|^\theta dx.\]

\begin{Lemma}
\label{L:Mit_EF}
Let $\theta>1$, $\alpha\in\R$ be given, such that $\mu_{\theta }(\alpha)>0$.
Let $q\ge \theta$ and assume that
$q\le \thetastar:=\frac{\theta n}{n-2\theta}$ if $n>2\theta$. Then there exists
a constant $c>0$ such that
$$
\irn|x|^\alpha|\Delta u|^\theta~\!dx\le 
c\left(\irn|x|^{-n+q\frac{n-2\theta+\alpha}{\theta}}|u|^q dx\right)^{\theta/q}
\quad\textit{for any $u\in \mathcal D^{2,\theta}(\R^n;|x|^\alpha dx)$.}
$$
\end{Lemma}

\proof
If $n>2\theta$ the conclusion readily follows via interpolation with the Sobolev inequality.
For a proof in lower dimensions, we use the Emden-Fowler transform
$T:C^2_c(\R^n\setminus\{0\})\to\in C^2_c(\R\times\S^{n-1})$, $T:u\mapsto g$ defined via
$$
u(x)=|x|^{\frac{2\theta-n-\alpha}{\theta }}~\!g\left(-\log|x|,\frac{x}{|x|}\right)~\!.
$$
We denote by
$\Delta_\sigma$  the Laplace-Beltrami operator on $\S^{n-1}$ and by
and $g'', g'$ the derivatives of $g=g(s,\sigma)$ with respect to $s\in\R$.
By direct computation one has that
\begin{gather*}
\irn |x|^{\alpha}|\Delta u|^\theta dx=\int_{\R}\int_{\S^{n-1}}
\left|\Delta_\sigma g+g''-2Ag'-\Gamma g\right|^\theta~dsd\sigma\\
\irn|x|^{\alpha-2\theta}| u|^\theta dx=\int_{\R}\int_{\S^{n-1}}
|g|^\theta~dsd\sigma~\!,
\end{gather*}
where
$\Gamma$ is given by (\ref{eq:gamma}) and
$ A=\frac{2(\theta -\alpha)+n(\theta -2)}{2\theta }$.
Thus, using the assumption $\mu_{\theta }(\alpha)>0$, one proves that
$$
\|g\|^\theta:=\int_{\R}\int_{\S^{n-1}}
\left|\Delta_\sigma g+g''-2Ag'-\Gamma g\right|^\theta~dsd\sigma
$$
is an equivalent norm on $W^{2,\theta}(\R\times\S^{n-1})$. Therefore,
$T$ can be regarded as an isometry between Banach spaces, and the conclusion 
readily follows by using the
 Sobolev embedding $W^{2,\theta}(\R\times\S^{n-1})\hookrightarrow
 L^q(\R\times\S^{n-1})$.
\QED

Under the assumptions in Lemma \ref{L:Mit_EF}, we have that the infimum
$$
S_{\theta ,q}(\alpha):  =  \inf_{\scriptstyle u\in \mathcal D^{2,\theta}(\R^n;|x|^\alpha dx)\atop\scriptstyle u\ne 0}\frac{\displaystyle
\int_{\R^n}|x|^{\alpha}|\Delta u|^\theta dx}
{\displaystyle\left(\int_{\R^n}|x|^{-n+q\frac{n-2\theta+\alpha}{\theta}}|u|^{q}dx\right)^{\theta /q}}~\!.
$$
is positive. Notice that for $n>2\theta$, $\alpha=0$ and $q=\theta^{*\!*}$ we have that
$$
S_{\theta ,\theta^{*\!*}}(0)=S^{*\!*}(\theta):  = 
 \inf_{\scriptstyle u\in \mathcal D^{2,\theta}(\R^n)\atop\scriptstyle u\ne 0}\frac{\displaystyle
\int_{\R^n}|\Delta u|^\theta dx}
{\displaystyle\left(\int_{\R^n}|u|^{\theta^{*\!*}}dx\right)^{\theta /\theta^{*\!*}}}~\!,
$$
which is the best constant in the Sobolev embedding
$\mathcal D^{2,\theta}(\R^n)\hookrightarrow L^{\theta^{*\!*}}(\R^n)$, see
\cite{Au}, \cite{Ta}.
The next existence results can be proved, for instance,
by using the techniques in \cite{CM1} (proof of Theorem 1.2). We omit details.

\begin{Theorem}
\label{T:general}
Let $\theta>1$, $\alpha\in\R$ be given, in such a way that
the infimum in (\ref{eq:Mit_ine}) is positive. Let $q>\theta$.

\begin{description}
\item$i)$
Assume that $n\ge 3$ and $q<\thetastar$ if $n>2\theta$.
Then $S_{\theta ,q}(\alpha)$ is achieved.

\item$ii)$
If $n>2\theta$ and
$S_{\theta ,\thetastar}(\alpha)<
S^{*\!*}(\theta)$
 then $S_{\theta ,\thetastar}(\alpha)$ is achieved.
\end{description}
\end{Theorem}

\begin{Remark}
\label{R:open}
Thanks to the results in \cite{Mit2} we know that 
$\mu_{\theta }(\alpha)=\left|\Gamma\right|^\theta>0$, whenever $\Gamma>0$.
We suspect that in this case the infimum $S_{\theta ,q}(\alpha)$ is always
achieved by radial functions.
We leave this as an open problem.
\end{Remark}

From Theorem \ref{T:general} one can easily infer sufficient conditions
for the existence of (minimal energy) solutions to the H\'enon-Lane-Emden
system (\ref{eq:HLE}), whenever $\mu_{p^\prime}(\alpha)>0.$
More can be said when $p=2$. From \cite{GM}, \cite{CM1} we know that
$\mu_{2}({-a})>0$ if and only if $-\Gamma$ is not an eigenvalue of the
Laplace-Beltrami operator on the sphere, where now
$$
\Gamma=\left(\frac{n+a}{2}\right)\left(\frac{n+b}{q}\right)~,\quad    \frac{a+n}{2}+\frac{b+n}{q}=n-2.
$$
From now on we assume that
\begin{equation}
\label{eq:non_resonance}
\textit{$-\left(\frac{n+a}{2}\right)\left(\frac{n+b}{q}\right)\neq k(n-2+k)$ for any integer $k\ge 0$}.
\end{equation}
By {\em ground state solutions} to 
(\ref{eq:p=2}) 
we mean
solutions $u,v$ such that $u$  achieves the
infimum
$$
\inf_{\scriptstyle u\in \mathcal D^{2,2}(\R^n;|x|^{-a}dx)\atop\scriptstyle u\ne 0}
\frac{\displaystyle
\int_{\R^n}|x|^{-a}\left|\Delta u\right|^{2} dx}
{\left(\displaystyle\int_{\R^n}|x|^b|u|^q~\!ds\right)^{2/q}}.
$$
For convenience of the reader we summarize here the main results
for (\ref{eq:p=2}) that can be obtained as immediate
corollaries of the results in \cite{CM2}.

\begin{Theorem}
\label{T:CM1}
Let $q>2$ and assume that $q<\starstar:=\frac{2n}{n-4}$ if $n\ge 5$.
\begin{description}
\item$~~i)$ If (\ref{eq:non_resonance}) holds, then
(\ref{eq:p=2}) has a ground state solution $\overline u, \overline v$.

\item$~ii)$ For every $q>2$ and for every integer
$k\ge 1$ there exists $\delta>0$ such that if
$$
0<|\Gamma+k(n-2+k)|<\delta
$$
then $\overline u$ is not radially symmetric. Thus problem (\ref{eq:HLE})
has at least two distinct weak solutions.

\item$iii)$ If
$|\Gamma|>\frac{n-1}{q-2}\left(1+\sqrt{q-1}\right)$
then $\overline u$ is not radially symmetric. Thus problem (\ref{eq:HLE})
has at least two distinct weak solutions.

\item$iv)$ Assume that
$-\Gamma>\frac{n-1}{2}$.
Then there exists $q_{\alpha}>2$ such that no ground state solution
to (\ref{eq:p=2}) can be positive.
\end{description}
\end{Theorem}
In the limiting case $n\ge 5$ and $q=\starstar$ the problem is more difficult.
We limit ourselves to point out some corollaries to the results in \cite{CM2}
in case $n\ge 6$.

\begin{Theorem}
Assume $n\ge 6$ and that (\ref{eq:non_resonance}) is satisfied. If in addition
$|a+2|>2$,
then  the problem
\begin{equation}
\label{eq:HLEstarstar}
\begin{cases}
-\Delta u=|x|^a v&\textrm{on $\R^n$}\\
-\Delta v=|x|^{b}|u|^{\frac{8}{n-4}}u&\textrm{on $\R^n$.}
\end{cases}
\end{equation}
has a ground state solution $\overline u, \overline v$. Moreover, the conclusions
$ii)$--$iv)$ in Theorem \ref{T:CM1} still hold. In particular, (\ref{eq:HLEstarstar}) has a radial and a non-radial weak solutions.
\end{Theorem}

\footnotesize

\label{References}

\end{document}